\documentclass[11pt]{amsart}
\usepackage{amsmath, amsfonts, amssymb, amstext, amscd, amsthm, graphicx, hyperref, url}
\usepackage[all,2cell,ps]{xy}
\allowdisplaybreaks

\newtheorem{goins-mugo}{Theorem}
\newtheorem{extension}[goins-mugo]{Corollary}
\newtheorem{elliptic_curve}[goins-mugo]{Proposition}
\newtheorem{3-points}[goins-mugo]{Corollary}
\newtheorem{elliptic_surface}[goins-mugo]{Proposition}
\newtheorem{4-points}[goins-mugo]{Theorem}

\begin{document}

\title{Points on Hyperbolas at Rational Distance}

\author{Edray Herber Goins}
\address{Purdue University \\ Department of Mathematics \\ 150 North University Street \\ West Lafayette, IN 47907-2067}
\email{egoins@math.purdue.edu}

\author{Kevin Mugo}
\address{Purdue University \\ Department of Mathematics \\ 150 North University Street \\ West Lafayette, IN 47907-2067}
\email{kmugo@math.purdue.edu}



\begin{abstract}
Richard Guy asked for the largest set of points which can be placed in the plane so that their pairwise distances are rational numbers.  In this article, we consider such a set of rational points restricted to a given hyperbola.

To be precise, for rational numbers $a$, $b$, $c$, and $d$ such that the quantity $D = \bigl( a \, d - b \, c \bigr) / \bigl( 2 \, a^2 \bigr)$ is defined and nonzero, we consider rational distance sets on the conic section $a \, x \, y + b \, x + c \, y + d = 0$.  We show that, if the elliptic curve $Y^2 = X^3 - D^2 \, X$ has infinitely many rational points, then there are infinitely many sets consisting of four rational points on the hyperbola such that their pairwise distances are rational numbers. We also show that any rational distance set of three such points can always be extended to a rational distance set of four such points.
\end{abstract}

\maketitle

\section{Introduction}

In Richard Guy's classic text \cite{Guy:2004p9886} there is a discussion of the largest set $S$ of points $P_t$ which can be placed in the plane $\mathbb A^2(\mathbb R)$ so that their pairwise distances $|| P_s - P_t ||$ are rational numbers.  One can find infinitely many such points if one restricts to a line.  Indeed, this is the case for the set $S = \{ \dots, \, P_t, \, \dots \}$ consisting of those points $P_t = \bigl( x_t : y_t : 1 \bigr)$ on the line $L: \, a \, x + b \, y = c$, all in terms of the rational numbers 
\[ \begin{aligned} x_t & = a \, c - b \, t & \qquad a & = \dfrac {2 \, r}{r^2 + 1} \\[2pt] y_t & = b \, c + a \, t & b & = \dfrac {r^2 - 1}{r^2 + 1} \end{aligned} \qquad \implies \qquad \left \| P_s - P_t \right \| = |s-t|. \]
\noindent  Similarly, one can find infinitely many such points if one restricts to a circle.  Indeed, this is the case for the set $S = \{ \dots, \, P_t, \, \dots \}$ consisting of those points $P_t = \bigl( x_t : y_t : 1 \bigr)$ on the circle $S^1: \, (x-h)^2 + (y-k)^2 = r^2$, all in terms of the rational numbers \[ \begin{aligned} x_t & = h + \frac {t^4 - 6 \, t^2 + 1}{(t^2 + 1)^2} \ r \\[2pt] y_t & = k + \frac {4 \, (t^3 - t)}{(t^2 + 1)^2} \ r \end{aligned} \quad \implies \quad \left \| P_s - P_t \right \| = \left| \frac {4 \, (s-t) \, (s \, t + 1)}{(s^2 + 1) \, (t^2 + 1)} \right| r. \]
\noindent Moreover, it can be shown that these rational distance sets are dense.

The problem becomes much more interesting if one restricts to those sets $S$ which are neither collinear nor concyclic.  Richard Guy asks in problem D20 of \cite{Guy:2004p9886} whether there exist six points $P_t$ in the plane $\mathbb A^2(\mathbb R)$, no three on a line $L$, no four on a circle $S^1$, all of whose mutual distances $|| P_s - P_t ||$ are rational; several authors have shown that infinitely many such sets $S$ exist.  In this article, we ask the same question, but we place the extra condition that these points are restricted to lie on a given conic section $\mathcal C$ -- so that, by B{\'e}zout's Theorem, no three will lie on a line and no five will lie on a circle.

Garikai Campbell \cite{MR2059753} showed that there exist infinitely many such sets $S$ of four points restricted to the parabola $C: \, y = x^2$.  He did so by showing that a certain elliptic curve has infinitely many rational points.  Motivated by this result and its proof, we show the following.

\begin{goins-mugo} For rational numbers $a$, $b$, $c$, and $d$ such that the quantity $D = \bigl( a \, d - b \, c \bigr) / \bigl( 2 \, a^2 \bigr)$ is defined and nonzero, consider the conic section
\[ \mathcal C: \qquad a \, x \, y + b \, x + c \, y + d = 0. \]
\noindent If the elliptic curve $E^{(D)}: \, Y^2 = X^3 - D^2 \, X$ has infinitely many rational points $(X: Y : 1)$, then there are infinitely many rational distance sets $S = \bigl \{ P_1, \, P_2, \, P_3, \, P_4 \bigr \}$ consisting of four rational points on $\mathcal C$ such that their pairwise distances $|| P_s - P_t ||$ are rational numbers. \end{goins-mugo}

\noindent In Theorem \ref{4-points} we show explicitly that we may choose the rational points
\[ \begin{aligned}
P_1 & = \left( - a \, c + (a \, d - b \, c) \, \dfrac {Y_1 \, X_2 \, X_3}{X_1 \, Y_2 \, Y_3} \ : \ - a \, b - a^2 \, \dfrac {X_1 \, Y_2 \, Y_3}{Y_1 \, X_2 \, X_3} \ : \ a^2 \right), \\[5pt]
P_2 & = \left( - a \, c + (a \, d - b \, c) \, \dfrac {X_1 \, Y_2 \, X_3}{Y_1 \, X_2 \, Y_3} \ : \ - a \, b - a^2 \, \dfrac {Y_1 \, X_2 \, Y_3}{X_1 \, Y_2 \, X_3} \ : \ a^2 \right), \\[5pt]
P_3 & = \left( - a \, c + (a \, d - b \, c) \, \dfrac {X_1 \, X_2 \, Y_3}{Y_1 \, Y_2 \, X_3} \ : \ - a \, b - a^2 \, \dfrac {Y_1 \, Y_2 \, X_3}{X_1 \, X_2 \, Y_3} \ : \ a^2 \right), \\[5pt]
P_4 & = \left( - a \, c + \dfrac {a \, d - b \, c}{4 \, D^2} \, \dfrac {Y_1 \, Y_2 \, Y_3}{X_1 \, X_2 \, X_3} \ : \ - a \, b - 4 \, D^2 \, a^2 \, \dfrac {X_1 \, X_2 \, X_3}{Y_1 \, Y_2 \, Y_3} \ : \ a^2 \right).
\end{aligned} \]

As a consequence of eliminating the $X_t$'s and $Y_t$'s from these formulas, we find that every rational distance set of three points can be extended to a rational distance set of four points:

\begin{extension} \label{extension} Say that $\{ P_1, \, P_2, \, P_3 \}$ is a rational distance set of three points $P_t = (x_t : y_t : 1)$ on the hyperbola $\mathcal C$.  If we choose the point
\[ P_4 = \begin{aligned} & \left( c + \dfrac {(a \, y_1 + b) \, (a \, y_2 + b) \, (a \, y_3 + b)}{a \, d - b \, c} \right. \\ & \left. \qquad \qquad \qquad : \ b + \dfrac {(a \, x_1 + c) \, (a \, x_2 + c) \, (a \, x_3 + c)}{a \, d - b \, c} \ : \ -a \right) \end{aligned} \]
\noindent then $S = \{ P_1, \, P_2, \, P_3, \, P_4 \}$ is a rational distance set of four points on the hyperbola.
\end{extension}

This work generalizes results found by Megan Ly, Shawn Tsosie, and Lyda Urresta in \cite{Ly:2010p19611} during the Mathematical Sciences Research Institute's Undergraduate Program.

\section{Main Results}

For the sequel, we fix rational numbers $a$, $b$, $c$, and $d$ such that the quantity $a \, (a \, d - b \, c) \neq 0$.  We also fix the notation
\[ \begin{aligned}
\mathcal C & = \left \{ (x:y:z) \in \mathbb P^2 \ \biggl| \ a \, x \, y + b \, x \, z + c \, y \, z + d \, z^2 =0 \right \}, \\[5pt]
D & = \dfrac {a \, d - b \, c}{2 \, a^2}, \\[5pt]
E^{(D)} & = \left \{ (X:Y:Z) \in \mathbb P^2 \ \biggl| \ Y^2 \, Z = X^3 - D^2 \, X \, Z^2 \right \}.
\end{aligned} \]

\noindent Let $S$ be rational distance set on $\mathcal C$, that is, let $S$ be a collection of affine rational points $P_t = (x_t : y_t : 1)$ on the hyperbola $\mathcal C$ such that their pairwise distances $|| P_i - P_j || = \sqrt{( x_i - x_j)^2 + ( y_i - y_j)^2}$ are rational numbers.   We give a way to generate such a set $S$.

\begin{elliptic_curve} \label{elliptic_curve} Fix an integer $n \geq 3$.  Consider a collection of $n \, (n-1)/2$ rational non-torsion points $(X_{ij} : Y_{ij} : 1)$ on the elliptic curve $E^{(D)}$ satisfying the $n \, (n-3)/2$ compatibility relations
\[ \dfrac {Y_{ij}}{X_{ij}} \,\dfrac {Y_{st}}{X_{st}} = \dfrac {Y_{is}}{X_{is}} \, \dfrac {Y_{jt}}{X_{jt}} \]
\noindent for all distinct indices $i$, $j$, $s$, and $t$.  Then the $n$ rational points
\[ P_t = \left( - a \, c + (a \, d - b \, c) \, \dfrac {Y_{ij} \, X_{it} \, X_{jt}}{X_{ij} \, Y_{it} \, Y_{jt}} \ : \ - a \, b - a^2 \, \dfrac {X_{ij} \, Y_{it} \, Y_{jt}}{Y_{ij} \, X_{it} \, X_{jt}} \ : \ a^2 \right) \]
\noindent for any distinct $i, \, j \neq t$ form a rational distance $S$ on the hyperbola $\mathcal C$. \end{elliptic_curve} 

\begin{proof} We explain how to derive this construction.  Let $S = \{ \dots, \, P_t, \, \dots \}$ be any rational distance set on $\mathcal C$.  For each pair $\{ P_i, \, P_j \}$ of distinct points in $S$, define the rational numbers
\[ X_{ij} = -\dfrac {D \, (x_i - x_j)}{( y_i - y_j) + \sqrt{( x_i - x_j)^2 + ( y_i - y_j)^2} }. \]
\noindent Upon solving for $y$ in the equation $a \, x \, y + b \, x + c \, y + d$, we find the slope
\[ \dfrac {X_{ij}^2 - D^2}{2 \, D \, X_{ij}} = \dfrac {y_i - y_j}{x_i - x_j} = \dfrac {a \, d - b \, c}{(a \, x_i + c) \, (a \, x_j + c)}. \]
\noindent If $|S| = n$, there are $\binom{n}{2} = n \, (n-1)/2$ such quadratic relations among the $x_t$, so we try and solve for the $\binom{n}{1} = n$ variables $x_t$ in terms of the $X_{ij}$.  We have the identity
\[ \dfrac {a \, d - b \, c}{(a \, x_t + c)^2} = \biggl[ \dfrac {X_{it}^2 - D^2}{2 \, D \, X_{it}} \biggr] \biggl[ \dfrac {X_{jt}^2 - D^2}{2 \, D \, X_{jt}} \biggr] \biggl[ \dfrac {X_{ij}^2 - D^2}{2 \, D \, X_{ij}} \biggr]^{-1} \quad \text{for distinct $i, \, j \neq t$.} \]
\noindent It makes sense to define the numbers
\[ n_{ij} = \dfrac {\sqrt{X_{ij}^3 - D^2 \, X_{ij}}}{X_{ij}} = \dfrac {a \, d - b \, c}{\sqrt{(a \, x_i + c) \, (a \, x_j + c)}} \]
\noindent which must satisfy the $\binom{n}{2} - \binom{n}{1} = n \, (n-3)/2$ compatibility conditions
\[ n_{ij} \, n_{st} = \dfrac {(a \, d - b \, c)^2}{\sqrt{(a \, x_i + c) \, (a \, x_j + c) \, (a \, x_s + c) \, (a \, x_t + c)}} = n_{is} \, n_{jt} \]
\noindent for all distinct indices $i$, $j$, $s$, and $t$.  Even though these numbers are \emph{not} rational in general, we still find the expression
\[ P_t = (x_t : y_t : 1) = \left( - a \, c + (a \, d - b \, c) \, \dfrac {n_{ij}}{n_{it} \, n_{jt}} \ : \ - a \, b - a^2 \, \dfrac {n_{it} \, n_{jt}}{n_{ij}} \ : \ a^2 \right) \]

\noindent for distinct $i, \, j \neq t$.

We now prove the result above.  Say that we are given a collection of $n \, (n-1)/2$ rational non-torsion points $(X_{ij} : Y_{ij} : 1)$ on the elliptic curve $E^{(D)}$ satisfying the $n \, (n-3) / 2$ compatibility relations above. Then the \emph{rational} numbers $n_{ij} = Y_{ij} / X_{ij}$ are nonzero, so that the $P_t$ as above are on the conic section.  Since we have the identity 
\[ \left \| P_s - P_t \right \| = \left| x_s - x_t \right| \sqrt{1 +\left( \dfrac {y_s - y_t}{x_s - x_t} \right)^2} = \biggl| \dfrac {n_{ij}}{n_{is} \, n_{js}} - \dfrac {n_{ij}}{n_{it} \, n_{jt}} \biggr| \, \biggl| \dfrac {X_{st}^2 + D^2}{X_{st}} \biggr| \]
\noindent we see that $S = \{ \dots, \, P_t, \, \dots \}$ is indeed a rational distance set.  \end{proof}

This result motivates study of the set
\[ E_n^{(D)} = \left \{ \prod_{i,j} \bigl( X_{ij} : Y_{ij} : Z_{ij}) \in \bigl( \mathbb P^2 \bigr)^{\frac {n \, (n-1)}{2}} \left| \begin{matrix} Y_{ij}^2 \, Z_{ij} = X_{ij}^3 - D^2 \, X_{ij} \, Z_{ij}^2 \\[5pt] X_{is} \, X_{jt} \, Y_{ij} \, Y_{st} = X_{ij} \, X_{st} \, Y_{is} \, Y_{jt} \end{matrix} \right. \right \} \] 
\noindent which is a projective variety of dimension $n \, (n-1)/2$.  The elliptic curve has involutions $E^{(D)} \to E^{(D)}$ defined by
\[ \left. \begin{aligned} \bigl( X' : Y' : Z' \bigr) & = [\pm 1] \, \bigl( X : Y : Z \bigr) \oplus \bigl( 0 : 0 : 1 \bigr) \\[5pt] & = \bigl( - D^2 \, X \, Z \, : \, \pm D^2 \, Y \, Z \, : \, X^2 \bigr) \end{aligned} \right \} \quad \implies \quad \dfrac {Y'}{X'} = \pm \dfrac {Y}{X}; \]

\noindent so the projective variety has involutions $E_n^{(D)} \to E_n^{(D)}$ as well.  This corresponds to the involution $n_{ij} \mapsto \pm n_{ij}$ in the formulas above.

We give some examples of this projective variety.

\textbf{Case $n = 3$.} Then $\binom{n}{2}  - \binom{n}{1} = 0$, so that $E_3^{(D)} = E^{(D)} \times E^{(D)} \times E^{(D)}$ is that projective variety of dimension 3 which is the product of elliptic curves.  
\[ \xymatrix @R=0.4in @C=0.5in {
& & E^{(D)} \\
E_3^{(D)} \ar@/_2pc/[drr] \ar@/^2pc/[rr] \ar@/^2pc/[urr] \ar@{->}[r] & E^{(D)} \times E^{(D)} \times E^{(D)} \ar@{->}[dr] \ar@{->}[r] \ar@{->}[ur] & E^{(D)} \\
& & E^{(D)} \\ 
& & (X_1 : Y_1 : Z_1) \\ {\left( \begin{matrix} (X_{23} : Y_{23} : Z_{23}), \\ (X_{13} : Y_{13} : Z_{13}), \\ (X_{12} : Y_{12} : Z_{12}) \end{matrix} \right)} \ar@/_2pc/[drr] \ar@/^2pc/[rr] \ar@/^2pc/[urr] \ar@{->}[r] & {\left( \begin{matrix} (X_1 : Y_1 : Z_1), \\ (X_2 : Y_2 : Z_2), \\ (X_3 : Y_3 : Z_3) \end{matrix} \right)} \ar@{->}[dr] \ar@{->}[r] \ar@{->}[ur] & (X_2 : Y_2 : Z_2) \\ & & (X_3 : Y_3 : Z_3) \\} \]

\noindent In other words, the formulae in Proposition \ref{elliptic_curve} uniquely determine $P_1$, $P_2$, and $P_3$.  We use this to generate rational distance sets $S$ on the hyperbola $\mathcal C$.

\begin{3-points} \label{3-points} If $(X_t : Y_t : 1)$ are three rational points on $E^{(D)}$ which are not points of finite order, then $S = \{ P_1, \, P_2, \, P_3 \}$ is a rational distance set on $\mathcal C: \ a \, x \, y + b \, x + c \, y + d = 0$ when we choose
\[ \begin{aligned}
P_1 = \left( - a \, c + (a \, d - b \, c) \, \dfrac {Y_1 \, X_2 \, X_3}{X_1 \, Y_2 \, Y_3} \ : \ - a \, b - a^2 \, \dfrac {X_1 \, Y_2 \, Y_3}{Y_1 \, X_2 \, X_3} \ : \ a^2 \right), \\[5pt]
P_2 = \left( - a \, c + (a \, d - b \, c) \, \dfrac {X_1 \, Y_2 \, X_3}{Y_1 \, X_2 \, Y_3} \ : \ - a \, b - a^2 \, \dfrac {Y_1 \, X_2 \, Y_3}{X_1 \, Y_2 \, X_3} \ : \ a^2 \right), \\[5pt]
P_3 = \left( - a \, c + (a \, d - b \, c) \, \dfrac {X_1 \, X_2 \, Y_3}{Y_1 \, Y_2 \, X_3} \ : \ - a \, b - a^2 \, \dfrac {Y_1 \, Y_2 \, X_3}{X_1 \, X_2 \, Y_3} \ : \ a^2 \right). 
\end{aligned} \]
\noindent In particular, when $E^{(D)}$ has positive rank, there are infinitely many rational distance sets $S$ of three points on the hyperbola $\mathcal C$. \end{3-points}

\begin{proof} The desired rational points
\[ \begin{aligned}
P_1 & = \left( - a \, c + (a \, d - b \, c) \, \dfrac {Y_{23} \, X_{12} \, X_{13}}{X_{23} \, Y_{12} \, Y_{13}} \ : \ - a \, b - a^2 \, \dfrac {X_{23} \, Y_{12} \, Y_{13}}{Y_{23} \, X_{12} \, X_{13}} \ : \ a^2 \right) \\[5pt] 
P_2 & = \left( - a \, c + (a \, d - b \, c) \, \dfrac {Y_{13} \, X_{12} \, X_{23}}{X_{13} \, Y_{12} \, Y_{23}} \ : \ - a \, b - a^2 \, \dfrac {X_{13} \, Y_{12} \, Y_{23}}{Y_{13} \, X_{12} \, X_{23}} \ : \ a^2 \right) \\[5pt] 
P_3 & = \left( - a \, c + (a \, d - b \, c) \, \dfrac {Y_{12} \, X_{13} \, X_{23}}{X_{12} \, Y_{13} \, Y_{23}} \ : \ - a \, b - a^2 \, \dfrac {X_{12} \, Y_{13} \, Y_{23}}{Y_{12} \, X_{13} \, X_{23}} \ : \ a^2 \right) \\[5pt] 
\end{aligned} \]
\noindent come about by relabeling the three points on $E^{(D)}$ as
\[ \begin{aligned}
\bigl( X_1 \, : \, Y_1 \, : \, Z_1 \bigr) & = \bigl( X_{23}: Y_{23}: Z_{23} \bigr) \\[5pt]
\bigl( X_2 \, : \, Y_2 \, : \, Z_2 \bigr) & = \bigl( X_{13}: Y_{13}: Z_{13} \bigr) \\[5pt]
\bigl( X_3 \, : \, Y_3 \, : \, Z_3 \bigr) & = \bigl( X_{12}: Y_{12}: Z_{12} \bigr) \\[5pt]
\end{aligned} \]
\noindent The result follows directly from Proposition \ref{elliptic_curve}. \end{proof}

\begin{figure}[t] \begin{center} \caption{Rational Distance Set on $x \, y + 12 = 0$} \label{fig:3-points} \includegraphics[width=0.95\textwidth]{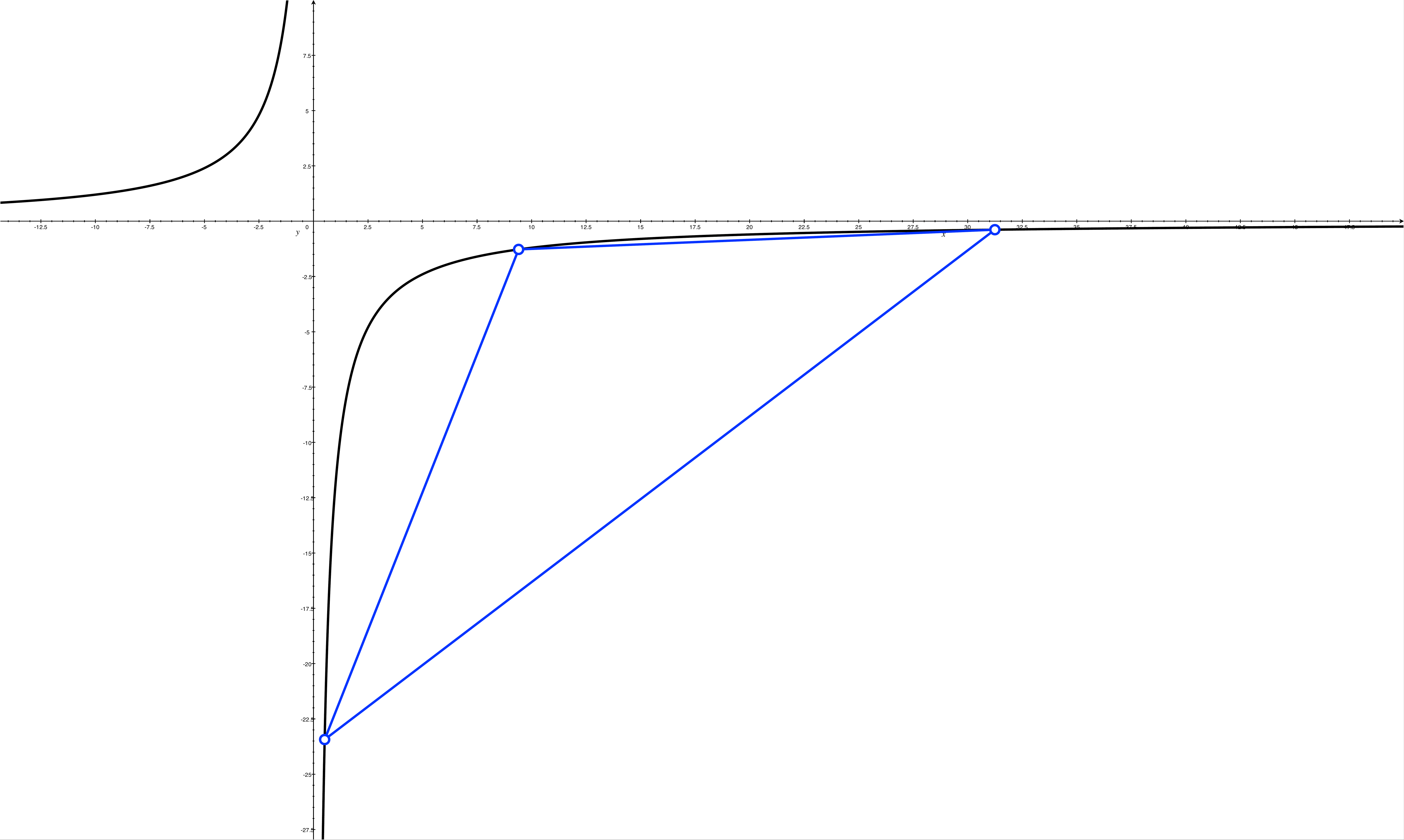} \end{center} \end{figure}

As a specific example, consider the conic section $\mathcal C: \, x \, y + 12 = 0$.  Then $D = 6$, and the elliptic curve $E^{(6)}$ has the three rational points $(12 : 36 : 1)$, $(50 : 35 : 8)$, and $(377844 : 2065932 : 12167)$; which correspond to $n_1 = 3$, $n_2 = 7/10$, and $n_3 = 4653/851$. respectively.  We find the affine rational points
\[ \left. \begin{aligned} P_1 & = \biggl( \dfrac {34040}{3619} \, : \, - \dfrac {10857}{8510} \, : \, 1 \biggr) \\[3pt] P_2 & = \biggl( \dfrac {11914}{23265} \, : \, - \dfrac {139590}{5957} \, : \, 1 \biggr) \\[3pt] P_3 & = \biggl( \dfrac {186120}{5957} \, : \, - \dfrac {5957}{15510} \, : \, 1 \biggr) \end{aligned} \right \} \quad \implies \quad \left \{ \begin{aligned} \left \| P_1 - P_2 \right \| & = \dfrac {1555297}{65142} \\[8pt] \left \| P_1 - P_3 \right \| & = \dfrac {28848020}{1319901} \\[8pt] \left \| P_2 - P_3 \right \| & = \dfrac {2129555051}{55435842} \end{aligned} \right. \]
\noindent lying on the hyperbola, so that $S = \{ P_1, \, P_2, \, P_3 \}$ is a rational distance set on $\mathcal C$. A plot of this configuration can be found in Figure \ref{fig:3-points}.

\textbf{Case $n = 4$.} Then $\binom{n}{2}  - \binom{n}{1} = 2$, so that $n_{12} \, n_{34} = n_{13} \, n_{24} = n_{14} \, n_{23}$ are the two compatibility conditions which define the projective variety $E_4^{(D)}$ of dimension 6.  We have a result that is similar to $E_3^{(D)} = E^{(D)} \times E^{(D)} \times E^{(D)}$.

\begin{elliptic_surface} \label{elliptic_surface} $E_4^{(D)}$ is the three-fold fiber product of the 2-dimensional variety
\[ \mathcal E^{(D)} = \left \{ \bigl( (u:v:w), \, T \bigr) \in \mathbb P^2 \times \mathbb P^1 \ \left| \ \begin{aligned} & v^4 \, w^3 = u \, \bigl( u + D^2 \, w \bigr) \\ & \times \bigl( u + D^2 \, T^4 \, w \bigr) \, \bigl( u^2 - D^4 \, T^4 \, w^2 \bigr)^2 \end{aligned} \right. \right \} \]
\noindent corresponding to the projection $\mathcal E^{(D)} \to \mathbb P^1$ which sends $\bigl( (u:v:w), \, T \bigr) \mapsto T$.
\end{elliptic_surface}

\begin{proof} We will explain the diagram
\[ \xymatrix @R=0.4in @C=0.5in {
& & \mathcal E^{(D)} \ar@{->}[dr] & \\
E_4^{(D)} \ar@/_2pc/[drr] \ar@/^2pc/[rr] \ar@/^2pc/[urr] \ar@{->}[r] & \mathcal E^{(D)} \times \mathcal E^{(D)} \times \mathcal E^{(D)} \ar@{->}[dr] \ar@{->}[r] \ar@{->}[ur] & \mathcal E^{(D)} \ar@{->}[r] & \mathbb P^1 \\
& & \mathcal E^{(D)} \ar@{->}[ur] & \\ 
& & \bigl( (u_1:v_1:w_1), \, T \bigr) \ar@{->}[dr] & \\ 
P \ar@/_2pc/[drr] \ar@/^2pc/[rr] \ar@/^2pc/[urr] \ar@{->}[r] & {\left( \begin{matrix} \bigl( (u_1:v_1:w_1), \, T \bigr), \\ \bigl( (u_2:v_2:w_2), \, T \bigr), \\ \bigl( (u_3:v_3:w_3), \, T \bigr) \end{matrix} \right)} \ar@{->}[dr] \ar@{->}[r] \ar@{->}[ur] & \bigl( (u_2:v_2:w_2), \, T \bigr) \ar@{->}[r] & T \\ 
& & \bigl( (u_3:v_3:w_3), \, T \bigr) \ar@{->}[ur] & \\ } \]

\noindent A rational point on the projective variety $E_4^{(D)}$ is in the form
\[ P = \left( \begin{aligned} & (X_1 : Y_1 : Z_1), \ (X_2 : Y_2: Z_2), \ (X_3:Y_3:Z_3), \\ & \qquad (X_1':Y_1':Z_1'), \ (X_2':Y_2':Z_2'), \  (X_3':Y_3':Z_3') \end{aligned} \right) \]
\noindent where $(X_t:Y_t:Z_t)$ and $(X_t': Y_t':Z_t')$ are six rational points on the elliptic curve $E^{(D)}$, that is,
\[ \begin{aligned} 
Y_1^2 \, Z_1 & = X_1^3 - D^2 \, X_1 \, Z_1^2 & \qquad {Y_1'}^2 \, Z_1' & = {X_1'}^3 - D^2 \, X_1' \, {Z_1'}^2 \\ 
Y_2^2 \, Z_2 & = X_2^3 - D^2 \, X_2 \, Z_2^2 & {Y_2'}^2 \, Z_2' & = {X_2'}^3 - D^2 \, X_2' \, {Z_2'}^2 \\
Y_3^2 \, Z_3 & = X_3^3 - D^2 \, X_3 \, Z_3^2 & {Y_3'}^2 \, Z_3' & = {X_3'}^3 - D^2 \, X_3' \, {Z_3'}^2 \\ 
\end{aligned} \]
\noindent and we have the two compatibility relations
\[ \dfrac {Y_1}{X_1} \, \dfrac {Y_1'}{X_1'} = \dfrac {Y_2}{X_2} \, \dfrac {Y_2'}{X_2'} = \dfrac {Y_3}{X_3} \, \dfrac {Y_3'}{X_3'}. \]

\noindent Given such a rational point $P$ on $E_4^{(D)}$, define the three rational points $\bigl( (u_t:v_t:w_t), \, T \bigr)$ on $\mathcal E^{(D)}$ in terms of
\[ \begin{aligned} 
u_t & = -4 \, D^3 \, X_t^2 \, X_t' \, Y_t' \, Z_t' \, \bigl(X_t^2 - D^2 \, Z_t^2 \bigr) \, \bigl({X_t'}^2 - D^2 \, {Z_t'}^2 \bigr) \\ 
v_t & = 8 \, D^5 \,X_t^2 \, {X_t'}^2 \, Z_t' \, \biggl[ D^2 \, Z_t^2 \, ( X_t' + D \, Z_t')^2 - X_t^2 \, (X_t' - D \, Z_t')^2 \biggr] \\
w_t & =  Y_t' \, \bigl(X_t^2 - D^2 \, Z_t^2 \bigr)^2 \, \bigl(X_t' - D \, Z_t' \bigr) \, \bigl(X_t' + D \, Z_t' \bigr)^3 \\[5pt]
T & = 2 \, D \, \dfrac {X_t}{Y_t} \, \dfrac {X_t'}{Y_t'}
\end{aligned} \]
\noindent where $T$ is independent of $t$.  Hence the expression
\[ Q = \bigl( (u_1: v_1: w_1), \ (u_2: v_2: w_2), \ (u_3: v_3: w_3), \ T \bigr) \]
\noindent represents a point on the three-fold fiber product of $\mathcal E^{(D)}$ with respect to the projection $\mathcal E^{(D)} \to \mathbb P^1$ which sends $\bigl( (u:v:w), \, T \bigr) \mapsto T$.

Conversely, given such a rational point $Q$ on the three-fold fiber product, define the rational points six rational points $(X_t:Y_t:Z_t)$ and $(X_t': Y_t':Z_t')$ on the elliptic curve $E^{(D)}$ in terms of
\begin{equation} \label{eq:elliptic_surface} \tag{$\ast$} \begin{aligned} 
X_t & = T \, v_t^2 \, w_t \\[5pt] 
Y_t & = v_t \, \bigl( u_t^2 - D^4 \, T^4 \, w_t^2 \bigr) \\[5pt] 
Z_t & = T^3 \, \bigl(u_t + D^2 \, w_t \bigr) \, \bigl( u_t^2 - D^4 \, T^4 \, w_t^2 \bigr) \\[5pt]
X_t' & = D \, w_t^2 \, \bigl( u_t^2 - D^4 \, T^4 \, w_t^2 \bigr) \\ & \quad \times \biggr[ 2 \, D \, v_t^2 \, w_t^2 - \bigl(u_t^2 - D^4 \, T^4 \, w_t^2 \bigr) \, \bigl(u_t^2 + 2 \, D^2 \, u_t \, w_t + D^4 \, T^4 \, w_t^2 \bigr) \biggl] \\[5pt]
Y_t' & = 2 \, D^2 \, v_t \, w_t^3 \\ & \quad \times \biggl[ 2 \, D \, v_t^2 \, w_t^2 - \bigl( u_t^2 - D^4 \, T^4 \, w_t^2 \bigr) \, \bigl( u_t^2 + 2 \, D^2 \, u_t \, w_t + D^4 \, T^4 \, w_t^2 \bigr) \biggr] \\[5pt]
Z_t' & = w_t^2 \, \bigl( u_t^2 - D^4 \, T^4 \, w_t^2 \bigr)^3
\end{aligned} \end{equation}

\noindent It is easy to verify that the two compatibility relations hold:
\[ \dfrac {Y_t}{X_t} = \dfrac {u_t^2 - D^4 \, T^4 \, w_t^2}{T \, v_t \, w_t}, \quad \dfrac {Y_t'}{X_t'} = \dfrac {2 \, D \, v_t \, w_t}{u_t^2 - D^4 \, T^4 \, w_t^2} \quad \implies \quad \dfrac {Y_t}{X_t} \, \dfrac {Y_t'}{X_t'} = \dfrac {2 \, D}{T} \]
\noindent so that we have a rational point $P$ on the projective variety $E_4^{(D)}$. \end{proof}

While the projective variety $E_4^{(D)}$ is somewhat cumbersome, the surface $\mathcal E^{(D)}$ gives us some insight on how to generate a family of rational distance sets $S$ consisting of four rational points $P_t$ on the hyperbola $\mathcal C$.

\begin{4-points} \label{4-points}  Continue notation as in Proposition \ref{elliptic_surface}.
\begin{enumerate}
\item There exists a nontrivial maps such that the composition
\[ \begin{CD} E^{(D)} @>>> \mathcal E^{(D)} @>>> \mathbb P^1 \end{CD} \]
\noindent is a constant map.
\item If $(X_t : Y_t : 1)$ are three rational points on $E^{(D)}$ which are not points of finite order, then $S = \{ P_1, \, P_2, \, P_3, \, P_4 \}$ is a rational distance set on $\mathcal C: \ a \, x \, y + b \, x + c \, y + d = 0$ when we choose
\[ \begin{aligned}
P_1 & = \left( - a \, c + (a \, d - b \, c) \, \dfrac {Y_1 \, X_2 \, X_3}{X_1 \, Y_2 \, Y_3} \ : \ - a \, b - a^2 \, \dfrac {X_1 \, Y_2 \, Y_3}{Y_1 \, X_2 \, X_3} \ : \ a^2 \right), \\[5pt]
P_2 & = \left( - a \, c + (a \, d - b \, c) \, \dfrac {X_1 \, Y_2 \, X_3}{Y_1 \, X_2 \, Y_3} \ : \ - a \, b - a^2 \, \dfrac {Y_1 \, X_2 \, Y_3}{X_1 \, Y_2 \, X_3} \ : \ a^2 \right), \\[5pt]
P_3 & = \left( - a \, c + (a \, d - b \, c) \, \dfrac {X_1 \, X_2 \, Y_3}{Y_1 \, Y_2 \, X_3} \ : \ - a \, b - a^2 \, \dfrac {Y_1 \, Y_2 \, X_3}{X_1 \, X_2 \, Y_3} \ : \ a^2 \right), \\[5pt]
P_4 & = \left( - a \, c + \dfrac {a \, d - b \, c}{4 \, D^2} \, \dfrac {Y_1 \, Y_2 \, Y_3}{X_1 \, X_2 \, X_3} \ : \ - a \, b - 4 \, D^2 \, a^2 \, \dfrac {X_1 \, X_2 \, X_3}{Y_1 \, Y_2 \, Y_3} \ : \ a^2 \right).
\end{aligned} \]
\item When $E^{(D)}$ has positive rank, there are infinitely many rational distance sets $S$ of four points on the hyperbola $\mathcal C$. 
\end{enumerate} \end{4-points}

\begin{proof} Define the map $E^{(D)} \to \mathcal E^{(D)}$ by
\[ (X: Y : Z) \mapsto \bigl( (X^2 \, Z \, : \, (X^2 + D^2 \, Z^2) \, Y \, : \, Z^3), \ 1 \bigr). \]
\noindent Clearly this is a nontrivial map whose composition with the projection $\mathcal E^{(D)} \to \mathbb P^1$ which sends $\bigl( (u:v:w), \, T \bigr) \mapsto T$ is a constant map.

Say that $(X_t : Y_t : 1)$ are three rational points on $E^{(D)}$ which are not points of finite order.  Then we have three points $\bigl( (u_t : v_t: z_t), \, 1 \bigr)$ on the surface $\mathcal E^{(D)}$ which corresponds to a rational point on the three-fold fiber product.  Using the expressions in equation \eqref{eq:elliptic_surface}, we find six rational points on $E^{(D)}$ given by
\[ \begin{aligned}
\bigl( X_{14}: Y_{14}: Z_{14} \bigr) & = \bigl( X_1' \, : \, Y_1' \, : \, Z_1' \bigr) \\ & = \bigl( D \, (X_1^2 - D^2 \, Z_1^2) \, : \, 2 \, D^2 \, Y_1 \, Z_1 \, : \, -(X_1 + D \, Z_1)^2 \bigr) \\
\bigl( X_{24}: Y_{24}: Z_{24} \bigr) & = \bigl( X_2' \, : \, Y_2' \, : \, Z_2' \bigr) \\ & = \bigl( D \, (X_2^2 - D^2 \, Z_2^2) \, : \, 2 \, D^2 \, Y_2 \, Z_2 \, : \, -(X_2 + D \, Z_2)^2 \bigr) \\
\bigl( X_{34}: Y_{34}: Z_{34} \bigr) & = \bigl( X_3' \, : \, Y_3' \, : \, Z_3' \bigr) \\ & = \bigl( D \, (X_1^2 - D^2 \, Z_1^2) \, : \, 2 \, D^2 \, Y_1 \, Z_1 \, : \, -(X_1 + D \, Z_1)^2 \bigr) \\
\end{aligned} \]
\noindent with the first three as labeled in Corollary \ref{3-points}.  We have chosen our maps so that $(Y_{ij}/X_{ij}) \, (Y_{st}/X_{st}) = 2 \, D$ for all distinct $i$, $j$, $s$, and $t$, so using the formulas in Proposition \ref{elliptic_curve} we have the four rational points as above, with the last being
\[ \begin{aligned}
P_4 & = \left( - a \, c + (a \, d - b \, c) \, \dfrac {Y_{12} \, X_{14} \, X_{24}}{X_{12} \, Y_{14} \, Y_{24}} \ : \ - a \, b - a^2 \, \dfrac {X_{12} \, Y_{14} \, Y_{24}}{Y_{12} \, X_{14} \, X_{24}} \ : \ a^2 \right) \\[5pt] 
& = \left( - a \, c + (a \, d - b \, c) \, \dfrac {Y_3 \, X_1' \, X_2'}{X_3 \, Y_1' \, Y_2'} \ : \ - a \, b - a^2 \, \dfrac {X_3 \, Y_1' \, Y_2'}{Y_3 \, X_1' \, X_2'} \ : \ a^2 \right) \\[5pt]
& = \left( - a \, c + \dfrac {a \, d - b \, c}{4 \, D^2} \, \dfrac {Y_1 \, Y_2 \, Y_3}{X_1 \, X_2 \, X_3} \ : \ - a \, b - 4 \, D^2 \, a^2 \, \dfrac {X_1 \, X_2 \, X_3}{Y_1 \, Y_2 \, Y_3} \ : \ a^2 \right) \\[5pt]
\end{aligned} \]
\noindent Using Corollary \ref{3-points}, we find precisely what is listed above.  \end{proof}

As a consequence, we find that we can always extend a rational distance set $S$ of three points on a hyperbola to a set of four points.

\vskip 0.1in

\noindent \textbf{Corollary 2.} \textit{Say that $\{ P_1, \, P_2, \, P_3 \}$ is a rational distance set of three points $P_t = (x_t : y_t : 1)$ on the hyperbola $\mathcal C: \ a \, x \, y + b \, x + c \, y + d = 0$.  If we choose the point
\[ P_4 = \begin{aligned} & \left( c + \dfrac {(a \, y_1 + b) \, (a \, y_2 + b) \, (a \, y_3 + b)}{a \, d - b \, c} \right. \\ & \left. \qquad \qquad \qquad : \ b + \dfrac {(a \, x_1 + c) \, (a \, x_2 + c) \, (a \, x_3 + c)}{a \, d - b \, c} \ : \ -a \right) \end{aligned} \]
\noindent then $S = \{ P_1, \, P_2, \, P_3, \, P_4 \}$ is a rational distance set of four points on the hyperbola.}

\begin{proof} Following the proof of Proposition \ref{elliptic_curve}, we can express the points $P_t$ in the form
\[ \begin{aligned}
P_1 = \left( - a \, c + (a \, d - b \, c) \, \dfrac {Y_1 \, X_2 \, X_3}{X_1 \, Y_2 \, Y_3} \ : \ - a \, b - a^2 \, \dfrac {X_1 \, Y_2 \, Y_3}{Y_1 \, X_2 \, X_3} \ : \ a^2 \right) \\[5pt]
P_2 = \left( - a \, c + (a \, d - b \, c) \, \dfrac {X_1 \, Y_2 \, X_3}{Y_1 \, X_2 \, Y_3} \ : \ - a \, b - a^2 \, \dfrac {Y_1 \, X_2 \, Y_3}{X_1 \, Y_2 \, X_3} \ : \ a^2 \right) \\[5pt]
P_3 = \left( - a \, c + (a \, d - b \, c) \, \dfrac {X_1 \, X_2 \, Y_3}{Y_1 \, Y_2 \, X_3} \ : \ - a \, b - a^2 \, \dfrac {Y_1 \, Y_2 \, X_3}{X_1 \, X_2 \, Y_3} \ : \ a^2 \right) 
\end{aligned} \]

\noindent where we have chosen the \emph{not necessarily rational} numbers
\[ \begin{aligned}
\dfrac {Y_1}{X_1} & = \dfrac {a \, d - b \, c}{a} \, \dfrac {1}{\sqrt{(a \, x_2 + c) \, (a \, x_3 + c)}} = \dfrac {\sqrt{(a \, y_2 + b) \, (a \, y_3 + b)}}{a} \\[5pt]
\dfrac {Y_2}{X_2} & = \dfrac {a \, d - b \, c}{a} \, \dfrac {1}{\sqrt{(a \, x_1 + c) \, (a \, x_3 + c)}} = \dfrac {\sqrt{(a \, y_1 + b) \, (a \, y_3 + b)}}{a}\\[5pt]
\dfrac {Y_3}{X_3} & = \dfrac {a \, d - b \, c}{a} \, \dfrac {1}{\sqrt{(a \, x_1 + c) \, (a \, x_2 + c)}} = \dfrac {\sqrt{(a \, y_1 + b) \, (a \, y_2 + b)}}{a}\\[5pt]
\end{aligned} \]

\noindent Following Theorem \ref{4-points}, we choose
\[ P_4 = \left( - a \, c + \dfrac {a \, d - b \, c}{4 \, D^2} \, \dfrac {Y_1 \, Y_2 \, Y_3}{X_1 \, X_2 \, X_3} \ : \ - a \, b - 4 \, D^2 \, a^2 \, \dfrac {X_1 \, X_2 \, X_3}{Y_1 \, Y_2 \, Y_3} \ : \ a^2 \right) \]
\noindent which is seen to be a rational point on $\mathcal C$.  \end{proof}

\begin{figure}[t] \begin{center} \caption{Rational Distance Set on $x \, y + 12 = 0$} \label{fig:4-points} \includegraphics[width=0.95\textwidth]{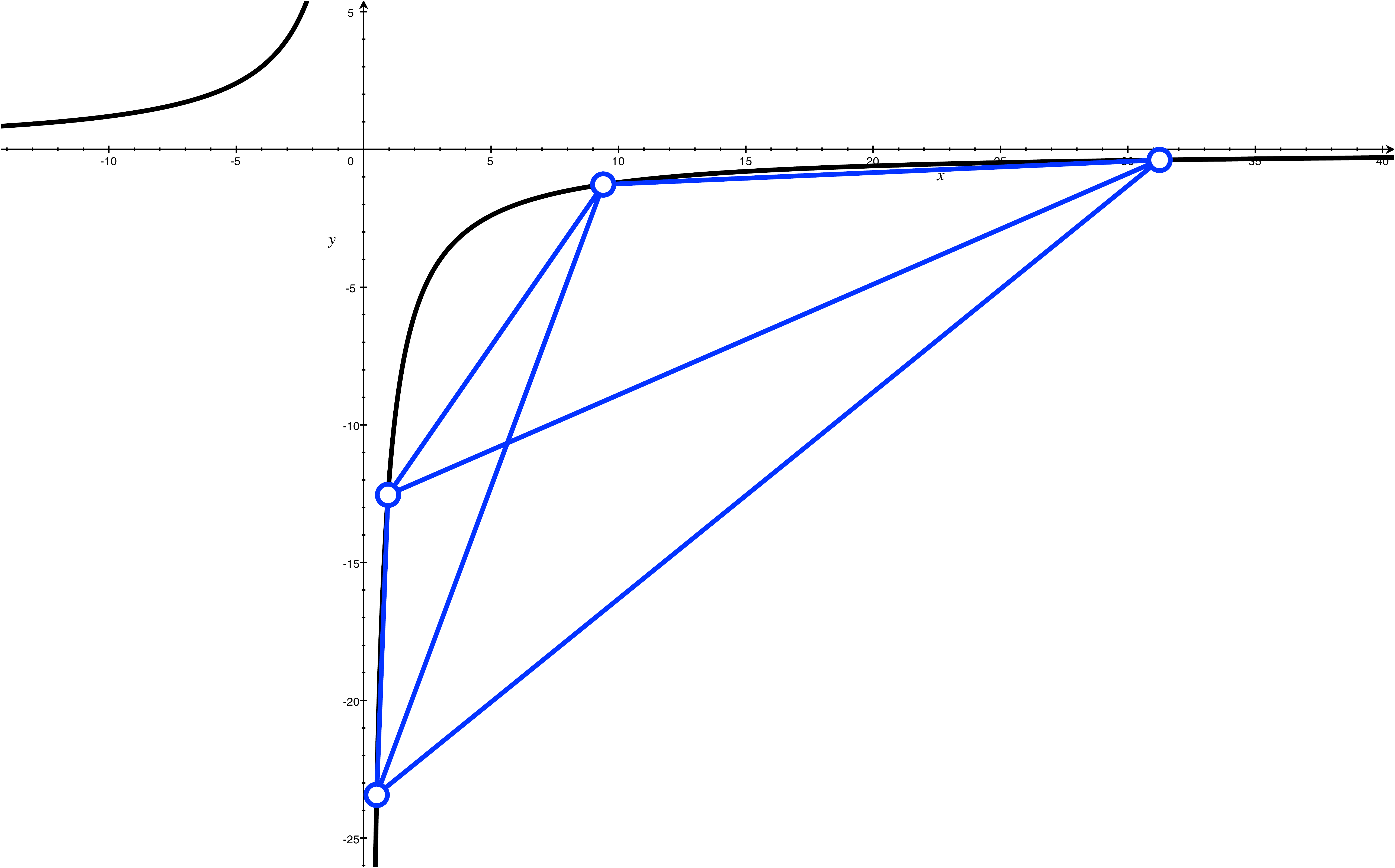} \end{center} \end{figure}

As a specific example, consider again the conic section $\mathcal C: \, x \, y + 12 = 0$.  We have seen that $D = 6$, and the elliptic curve $E^{(6)}$ has the three rational points $(12 : 36 : 1)$, $(50 : 35 : 8)$, and $(377844 : 2065932 : 12167)$.  These correspond to the rational points $\bigl( (144:6480:1), \, 1 \bigr)$, $\bigl( (5000:42035:128), \, 1 \bigr)$, and $\bigl( (3283620031728:578364811524720:3404825447), \, 1 \bigr)$, respectively, on the surface $\mathcal E^{(D)}$.  We find the affine rational points
\[ \left. \begin{aligned} P_1 & = \biggl( \dfrac {34040}{3619} \, : \, - \dfrac {10857}{8510} \, : \, 1 \biggr) \\[3pt] P_2 & = \biggl( \dfrac {11914}{23265} \, : \, - \dfrac {139590}{5957} \, : \, 1 \biggr) \\[3pt] P_3 & = \biggl( \dfrac {186120}{5957} \, : \, - \dfrac {5957}{15510} \, : \, 1 \biggr) \\[5pt] P_4 & = \biggl( \dfrac {32571}{34040} \, : \, - \dfrac {136160}{10857} \, : \, 1 \biggr) \end{aligned} \right \} \quad \implies \quad \left \{ \begin{aligned} \left \| P_1 - P_2 \right \| & = \dfrac {1555297}{65142} \\[8pt] \left \| P_1 - P_3 \right \| & = \dfrac {28848020}{1319901} \\[8pt] \left \| P_2 - P_3 \right \| & = \dfrac {2129555051}{55435842} \\[5pt]  \left \| P_1 - P_4 \right \| & = \dfrac {1040847151}{73914456} \\[5pt] \left \| P_2 - P_4 \right \| & = \dfrac {1726556399}{158388120} \\[5pt] \left \| P_3 - P_4 \right \| & = \dfrac {1555297}{47656} \end{aligned} \right. \]
\noindent lying on the hyperbola, so that $S = \{ P_1, \, P_2, \, P_3, \, P_4 \}$ is a rational distance set on $\mathcal C$.  A plot of this configuration can be found in Figure \ref{fig:4-points}.

\end{document}